\documentclass[12pt,a4paper]{amsart}
\usepackage[top=35mm, bottom=35mm, left=30mm, right=30mm]{geometry}
\usepackage[colorlinks=true,citecolor=blue]{hyperref}
\usepackage{mathptmx}
\usepackage{eucal}
\usepackage{graphicx}
\usepackage{mathrsfs}
\usepackage{amssymb}
\usepackage{amsmath}
\usepackage{amsthm}
\usepackage{amscd}
\usepackage{xcolor}
\usepackage{verbatim}

\newtheorem{thm}{Theorem}[section]

\newtheorem{lem}[thm]{Lemma}

\theoremstyle{definition}

\numberwithin{equation}{section}

\newcommand{\N}{\mathbb{N}}

\def \M { \mathcal{M}}

\numberwithin{equation}{section}

\input xy
\xyoption{all}

\begin{document}
	\title[Discrete spectrum of probability measures]{Discrete spectrum of probability measures for  locally compact  group actions}
	\author[Z. Hu, X. Ma, L. Xu and X. Zhou]
	{Zongrui Hu, Xiao Ma, Leiye Xu and Xiaomin Zhou}
	
	\address[Z. Hu, X. Ma and L. Xu]
	{CAS Wu Wen-Tsun Key Laboratory of Mathematics, School of Mathematical Sciences, University of Science and
		Technology of China, Hefei, Anhui 230026, China}
	\email{zongrui@mail.ustc.edu.cn, xiaoma@ustc.edu.cn, leoasa@mail.ustc.edu.cn}

	\address[X.~Zhou]{School of Mathematics and Statistics, Huazhong University of Science and Technology, Wuhan, Hubei 430074, China}
	\email{zxm12@mail.ustc.edu.cn}
	
    \thanks{Corresponding author: Xiaomin Zhou.}
	
	

\begin{abstract}In this paper, we investigate the discrete spectrum of probability measures for actions of locally compact groups. We establish that a probability measure has a discrete spectrum if and only if it has bounded measure-max-mean-complexity. 

As applications:
1) An invariant measure for a locally compact amenable group action has a discrete spectrum if and only if it has bounded mean-complexity along F\o lner sequences;
2) An invariant measure for a locally compact amenable group action has a discrete spectrum if and only if it is mean equicontinuous along a tempered F\o lner sequence, or equicontinuous in the mean along a tempered F\o lner sequence.
\end{abstract}

\maketitle

	\section{Introduction}
	Throughout this paper,  by a \textbf{locally compact group action} we mean a pair $(X,G)$ where $X$ is a compact metric space with  metric $d$ and $G$ is a locally compact group acting on $X$ continuously. Let $\M(X)$ denote the collection of all Borel probability measures on a measurable space $(X,\mathcal{B}_X)$, where $\mathcal{B}_X$ denotes the Borel $\sigma$-algebra on $X$. If $(X,G)$ is a locally compact group action,  the set of all $G$-invariant Borel probability measures in $\M(X)$ will be denoted by $\M(X,G)$.  By a measure-preserving system  we mean a quadruple $(X, \mathcal{B}_X, \mu, G)$, where $(X,G)$ is a locally compact  group action and $\mu$ is a $G$-invariant Borel probability measure on $X$.

	  For any invertible measure-preserving system $(X, \mathcal{B}_{X}, \mu, T)$, the space of $L^2$ functionals, i.e., $L^2(X,\mathcal{B}_X, \mu)$, can be decomposed into two subspaces: one consisting of functionals with continuous spectrum and the other with discrete spectrum.  This decomposition is a consequence of the Koopman-von Neumann spectrum mixing theorem \cite{KV32}. For an infinite discrete countable group action $(X, G)$ and its corresponding $G$-invariant Borel probability measures, a similar result was established by Bergelson in \cite{Be}. By a system has discrete spectrum, we mean that there exists an orthonormal basis of eigenfunctions for the space $L^2(X,\mathcal{B}_X, \mu)$. In a sense, such systems are the simplest possible dynamical systems.
	  
There is a long history for the characterizations of discrete spectrum, which has attracted substantial interest in the literature since its introduction in the field of dynamical systems. 
For example, in 1967, Ku\v{s}nirenko \cite{Ku67} introduced the concept of measure-theoretical sequence entropy for $\mathbb{Z}$-actions and demonstrated that a $\mathbb{Z}$-action has discrete spectrum if and only if its sequence entropy is zero with respect to any infinite subsequence of $\mathbb{N}\cup\{0\}$.
	In 1982, Scarpellini \cite{S} revealed that a large class of measure preserving flows with  discrete spectrum exhibit some Lyapunov-like stability properties.
	The relationship between discrete spectrum and pure point diffraction has also been extensively explored. Lenz and Strungaru \cite{LS} established the equivalence between discrete spectrum and pure point diffraction for measure dynamical systems on locally compact $\sigma$-compact abelian groups. 
This result was subsequently extended to actions of locally compact $\sigma$-compact abelian groups on compact metric spaces by Lenz \cite{L} and by Lenz, Spindeler, and Strungaru \cite{LSS}. Recently, Fuhrmann, Gröger, and Lenz \cite{FGL} demonstrated this equivalence for actions of locally compact $\sigma$-compact amenable groups on compact metric spaces.
	
	For certain simple dynamical systems—such as those with zero entropy or those possessing discrete spectrum—the complexity function emerges as a valuable tool for analysis.
	The investigation of complexity functions dates back to the pioneering work of Morse and Hedlund \cite{HM}, who explored the complexity functions of symbolic systems and demonstrated that the boundedness of these functions is equivalent to the eventual periodicity of the system. 
	Katok \cite{K} subsequently introduced the complexity function through a modified notion of spanning sets, incorporating an invariant measure and an error parameter. 
	Building on this foundation, Huang, Wang, and Ye \cite{HWY} defined the measure complexity of an invariant measure using the mean metric. They established that if an invariant measure of a transformation has  discrete spectrum, its measure complexity must be bounded. 
	This result was later complemented by Huang, Li, Thouvenot, Xu, and Ye \cite{HLTXY}, who proved the converse: bounded measure complexity implies discrete spectrum.
Further advancements were made by Ferenczi \cite{Fe} and Yu \cite{Yu}, who employed the $\alpha$-names of partitions and the Hamming distance to show that the complexity function of a system is bounded if and only if the system has discrete spectrum. 
This characterization was extended to countable discrete amenable group actions by Yu, Zhang, and Zhang \cite{YZZ}, who introduced measure complexity along F\o lner sequences and demonstrated that an invariant measure has discrete spectrum if and only if its measure complexity is bounded along any F\o lner sequence. Xu and Xu \cite{XX} further generalized this result, proving that for general group actions, an invariant measure has discrete spectrum if and only if it has finite max-mean-measure complexity.
These developments naturally lead to the following question: 
\smallskip

\noindent{\bf Question:} {\it Can we characterize discrete spectrum using measure complexity in the context of group actions involving locally compact  groups over compact metric spaces, following the approaches of \cite{HLTXY, HWY, YZZ}?}
\smallskip

\noindent In this paper, we provide an affirmative answer, affirming the utility of measure complexity as a robust tool for characterizing discrete spectrum in broader dynamical settings.

With the above question in mind, we consider Borel probability measures for locally compact group actions, and investigate the connection between discrete spectrum and measure complexity.
More precisely,  for a locally compact group action $(X,G)$, we say $\mu\in\M(X)$, which is a Borel probability measure on $X$,  has \emph{discrete spectrum} if for any $f\in C(X)$, $\{gf:g\in G\}$ is precompact in $L^2(\mu)$, where $(gf)(x):=f(gx)$ for any $g\in G$ and all $x\in X.$ 
Here we use continuous functions to define the discrete spectrum mainly to ensure the measurability after transformations. 
This definition is equivalent to the usual definition when $(X,G)$ is a topological dynamical system and $\mu\in\M(X,G)$.

Let $(X,G)$ be a  locally compact group action and $\mu\in\mathcal{M}(X)$.
For $\rho\in\M(G)$,  where $\mathcal{M}(G)$ is the space  of Borel probability measures on $G$, and $x,y\in X$, define
$$ d_\rho(x,y)=\int_{G}d(gx,gy)d\rho(g).$$
For $\epsilon>0$ and $x\in X$, define
$$ B_\rho(x,\epsilon)=\{y\in X: d_\rho(x,y)<\epsilon\}.$$
Let $S_\rho(X,\mu,G,\epsilon)$ denote the minimal natural number $K$ such that there exists a subset  
$$\{x_1,x_2,\dots,x_K\}\subseteq X$$ such that
$$\mu(\cup_{k=1}^K  B_\rho(x_k,\epsilon))>1-\epsilon.$$
We denote the collection of all finite subsets of $G$ by $F(G)$. Denote $\rho_E=\frac{1}{|E|}\sum_{g\in E}\delta_g$ for $E\in F(G)$, where $\delta_g$ is the Dirac measure.

We say $\mu$ has \emph{bounded measure-max-mean-complexity} if for any $\epsilon>0$ there is $K(\epsilon)>0$ such that ${S}_\rho(X, \mu, G, \epsilon)<K(\epsilon)$ for every  $\rho\in \M(G)$. 	We say $\mu$ has \emph{bounded max-mean-complexity}, if for any $\epsilon>0$ there is $K(\epsilon)>0$ such that ${S}_{\rho_E}(X, \mu, G, \epsilon)<K(\epsilon)$ for every  $E\in F(G)$.

For locally compact group actions and Borel probability measures, we obtain the following equivalent characterizations of discrete spectrum.
\begin{thm}\label{thm-2} Let $(X,G)$ be a  locally compact group action and $\mu\in \M(X)$. Then the followings are equivalent.
	\begin{itemize}
		\item[(1)] $\mu$ has  discrete spectrum;
		\item[(2)] $\mu$ has bounded measure-max-mean-complexity;
		\item[(3)] $\mu$ has bounded max-mean-complexity.
	\end{itemize}
\end{thm}

As an application of Theorem \ref{thm-2}, we also obtain the equivalent characterizations of discrete spectrum
for locally compact amenable group actions.
Recall that a locally compact group $G$ is called amenable (see \cite{F}) if
there exists a sequence of non-empty compact sets $\{F_n\}_{n\in\N}$ in $G$  with $0< m_G(F_n)<\infty$ for all $n\in \mathbb{N}$ such that
\begin{align*}
	\lim_{n\to\infty}\frac{m_G(gF_n\triangle F_n)}{m_G(F_n)}
	=0\quad\textnormal{for all }g\in G,
\end{align*}
where $\triangle$ denotes the symmetric difference and $m_G$ is a
\emph{(left) Haar measure} of $G$. And the sequence $\{F_n\}_{n\in \mathbb{N}}$ of compact subsets of $G$ is called  a \emph{(left) F\o lner sequence} of $G$.

To articulate the subsequent result with exactness, we require some preliminary work.
Define $m_G|_{F_n}(\cdot):=\frac{m_G(\cdot\cap F_n)}{m_G( F_n)}.$
We say $\mu$ has \emph{bounded mean-complexity along a F\o lner sequence $\{F_n\}_{n\in\N}$  of $G$}, if for any $\epsilon>0$ there is $K(\epsilon)>0$, such that for all $n\in \mathbb{N}$, one has $$S_{m_G|_{F_n}}(X,\mu,G,\epsilon)<K(\epsilon),$$
where $S_{m_{G}|_{F_n}}(X, \mu, G, \epsilon)$ is the minimal natural number $K$ satisfying that there exists a subset $$\{x_{1},\dots, x_K\}\subseteq X$$ such that $$\mu(\cup_{k=1}^{K}B_{m_G|_{F_n}}(x_k, \epsilon):=\{x\in X: d_{m_{G}|_{F_n}}(x,x_k)<\epsilon\})>1-\epsilon,$$
 and $$d_{m_{G}|_{F_n}}(x,y):= \int_{G}d(gx, gy) dm_G|_{F_n}(g)=\frac{1}{m_G(F_n)}\int_{F_n}d(gx, gy)dm_G(g) \text{ for any } x,y\in X.$$ 
\begin{thm}\label{thm-1}Let $G$ be a locally compact amenable group and $(X,\mathcal{B}_X,\mu,G)$ be a measure-preserving system. Then the followings are equivalent.
	\begin{itemize}
		\item[(1)] $\mu$ has bounded mean-complexity along some F\o lner sequence of $G$;
		\item[(2)] $\mu$ has bounded mean-complexity along any F\o lner sequence of $G$;
		\item[(3)] $\mu$ has bounded measure-max-mean-complexity;
		\item[(4)] $\mu$ has discrete spectrum.
	\end{itemize}
\end{thm}

As a further application of  Theorem \ref{thm-2}, we establish the equivalence between discrete spectrum and mean equicontinuous along a tempered F\o lner sequence (equicontinuous in the mean along a tempered F\o lner sequence) for a measure-preserving system $(X, \mathcal{B}_X, G,\mu)$, where $G$ is a locally compact amenable group. 
This result affirmatively addresses Question 6 posed by Li, Ye, and Yu in their comprehensive survey \cite{LYY}.

A F\o lner sequence $\{F_n\}_{n\in \mathbb{N}}$ in $G$ is said to be tempered (see Shulman \cite{Sh}) if there exists a positive constant $C$ which is independent of $n$ such that for any $n\in \mathbb{N}$
\begin{align}\label{sc}
|\bigcup_{k<n}F_k^{-1}F_n|\leq C|F_n|.
\end{align}
In Linderstrauss \cite{Li}, (\ref{sc}) is also called Shulman Condition.
Lindenstrauss in \cite{Li} developed general covering lemmas and using these tools to prove the
pointwise ergodic theorem for general amenable groups.

We say that a subset $Q\subseteq X$ is \emph{mean equicontinuous along $\{F_n\}_{n\in  \mathbb{N}}$ (equicontinuous in the
mean along $\{F_n\}_{n\in  \mathbb{N}}$)} if for any $\epsilon>0$, there exists $\delta>0$ such that
 $$x,y\in Q \text{ with } d(x,y)<\delta \Rightarrow \limsup_{n\to\infty}d_{m_G|_{F_n}}(x,y)<\epsilon\ (d_{m_G|_{F_n}}(x,y)<\epsilon, \forall n\in \mathbb{N}).$$
We say $\mu$ is \emph{mean equicontinuous along $\{F_n\}_{n\in  \mathbb{N}}$ (equicontinuous in the
mean  along $\{F_n\}_{n\in  \mathbb{N}}$)} if for any $\tau>0$, there exists a compact
subset $Q$ of $X$ with $\mu(Q)>1-\tau$ such that $Q$ is mean equicontinuous along $\{F_n\}_{n\in  \mathbb{N}}$ (equicontinuous in the
mean along $\{F_n\}_{n\in  \mathbb{N}}$) .

\begin{thm}\label{thm-3}Let $G$ be a locally compact amenable group  and $(X,\mathcal{B}_X, \mu,G)$ be a measure-preserving system with a tempered F{\o}lner sequence $\{F_n\}_{n\in \mathbb{N}}.$   Then the followings are equivalent:
		\begin{itemize}
			\item[(1)] $\mu$ has discrete spectrum;
			\item[(2)] $\mu$ is mean equicontinuous along $\{F_n\}_{n\in \mathbb{N}};$

             \item[(3)] $\mu$ is equicontinuous in the mean along $\{F_n\}_{n\in  \mathbb{N}}.$
		\end{itemize}
	\end{thm}

The outline of the paper is the following. In Section \ref{s-t-1}, we give the proof of Theorem \ref{thm-2}. In Section \ref{s-t-2},  we present the proof of Theorem \ref{thm-1}. In Section \ref{s-t-3}, the proof of Theorem \ref{thm-3} is provided.

	\section{Proof of Theorem \ref{thm-2}}\label{s-t-1}
	 The aim of this section is to prove Theorem \ref{thm-2}, i.e., the equivalence between the discrete spectrum, the bounded measure-max-mean-complexity and the bounded max-mean-complexity of a Borel probability measure $\mu\in \mathcal{M}(X)$ for a locally compact group action $(X, G)$.
	Before proving Theorem \ref{thm-2}, we need the following lemma.  

For a Borel probability measure $\mu\in \M(X)$, we say the product measure $\mu\times \mu\in \mathcal{M}(X\times X)$ has discrete spectrum, if for any $F\in C(X\times X)$, $\{gF: g\in G\}$ is precompact, where $gF(x,y)=F(gx, gy)$ for any $x,y\in X$.
	\begin{lem}\label{l-1}Let $(X,G)$ be a  locally compact group action and $\mu\in \M(X)$. If $\mu$ has discrete spectrum, so is $\mu\times\mu\in \mathcal{M}(X\times X)$.
	\end{lem}
\begin{proof}Given $F\in C(X\times X)$ and $\epsilon>0$. There exists $K\in\N$ and $ f_1,\dots,f_K,\tilde f_1,\dots,\tilde f_K\in C(X)$ such that $$\|F-\sum_{k=1}^Kf_k\otimes \tilde f_k\|_\infty<\frac{\epsilon}{2},$$
where $f\otimes \tilde f(x,y):=f(x)\cdot\tilde f(y)$ for any $(x,y)\in X\times X.$
	Since $\mu$ has discrete spectrum, there exist finite subsets $H_k,\tilde H_k$ of $G$ for $k=1,2,\dots,K$ such that for every $g\in G$, there exists $h_k^{(g)}\in H_k,\tilde h_k^{(g)}\in \tilde H_k$ such that
	$$\|gf_k-h_k^{(g)}f_k\|_{L^2(\mu)}\le \frac{\epsilon}{4K(\|\tilde f_k\|_\infty+\|f_k\|_\infty)}$$ and $$ \|g\tilde f_k-\tilde h_k^{(g)}\tilde f_k \|_{L^2(\mu)}\le \frac{\epsilon}{4K(\|\tilde f_k\|_\infty+\|f_k\|_\infty)}$$ for $k=1,2,\dots,K$.
Thus for $g\in G$,
\begin{align*}
&\|gF-\sum_{k=1}^Kh_k^{(g)}f_k\otimes\tilde h_k^{(g)} \tilde f_k\|_{L^2(\mu\times\mu)}\\
&\le \|gF-\sum_{k=1}^Kgf_k\otimes g\tilde f_k\|_\infty+\|\sum_{k=1}^Kgf_k\otimes g\tilde f_k-\sum_{k=1}^Kh_k^{(g)}f_k\otimes\tilde h_k^{(g)} \tilde f_k\|_{L^2(\mu\times\mu)}\\
&< \frac{\epsilon}{2}+\sum_{k=1}^K\|gf_k\otimes g\tilde f_k-h_k^{(g)}f_k\otimes\tilde h_k^{(g)} \tilde f_k\|_{L^2(\mu\times\mu)}\\
&\le \frac{\epsilon}{2}+\sum_{k=1}^K\left(\|gf_k\otimes g\tilde f_k-h_k^{(g)}f_k\otimes g\tilde f_k\|_{L^2(\mu\times\mu)}+\|h_k^{(g)}f_k\otimes g\tilde f_k-h_k^{(g)}f_k\otimes \tilde h_k^{(g)} \tilde f_k\|_{L^2(\mu\times\mu)}\right)\\
&\le  \frac{\epsilon}{2}+\sum_{k=1}^K\left(\|\tilde f_k\|_\infty\|gf_k-h_k^{(g)}f_k\|_{L^2(\mu)}+\| f_k\|_\infty\| g\tilde f_k-\tilde h_k^{(g)} \tilde f_k\|_{L^2(\mu)}\right)\\
&\le \epsilon.
\end{align*}
 	Denote $\mathcal{N}_\epsilon=\{\sum_{k=1}^Kh_kf_k\otimes\tilde h_k \tilde f_k:\forall h_k\in H_k,\forall\tilde h_k\in\tilde H_k, k=1,\dots,K \}$.  Then $\mathcal{N}_\epsilon$ is an $\epsilon$-net of $\{gF:g\in G\}$.
 	By the arbitrariness of $\epsilon$, we can deduce that $\{gF:g\in G\}$ is precompact in $L^2(\mu\times\mu)$. Therefore, $\mu\times\mu$ has discrete spectrum.	
\end{proof}
Now we are ready to prove Theorem \ref{thm-2}.
 \begin{proof}[Proof of Theorem \ref{thm-2}]There are three parts of the proof, i.e.,
\begin{itemize}
\item[(1)$\Rightarrow$(2):]  $\mu\in\mathcal{M}(X)$ has bounded measure-max-mean-complexity if it has discrete spectrum;
\item[(2)$\Rightarrow$(3):] $\mu\in \mathcal{M}(X)$ with bounded measure-max-mean-comlexity also has bounded max-mean-complexity; 
\item[(3)$\Rightarrow$(1):] if $\mu\in \mathcal{M}(X)$ has bounded max-mean-complexity, then it has discrete spectrum.
\end{itemize}
Note that (2)$\Rightarrow$(3) holds straightforwardly, according to the definition of bounded measure-max-mean-complexity and bounded max-mean-complexity. Therefore, we only prove (1)$\Rightarrow$(2) and (3)$\Rightarrow$(1) in the following.
\smallskip

\noindent{\bf (1)$\Rightarrow$(2):} Without loss of generality, we can assume that \begin{align}\label{diam1}\text{diam}(X)=1.\end{align} Assume that $\mu\in\mathcal{M}(X)$  has discrete spectrum. Then so is $\mu\times\mu\in \mathcal{M}(X\times X)$ by Lemma \ref{l-1}. Put $f(x,y)=d(x,y)$ on $X\times X$. Fix $\delta>0$. 	We can find $K\in\mathbb{N}$ and  a finite Borel measurable partition $\beta=\{B_1,B_2,\dots,B_K\}$ of $X$ with  $\text{diam}(B_k)<\delta$ for $ 1\leq k\leq K$. Since $\mu\times \mu$ has discrete spectrum, there exists a finite subset $H$ of $G$ and a Borel measurable map $\varphi: G\to H$ such that
 \begin{align}\label{eq-12}\|gf-\varphi(g)f\|_{L^2(\mu\times\mu)}\le \frac{\delta^{4}}{K}.
 \end{align}	
We are going to show that $$ S_\rho(X,\mu,G,6\delta+4\sqrt{\delta})\le K^{|H|}$$ for every $\rho\in\M(G)$ in three steps. 
\smallskip

\noindent{\bf Step 1:} At first, we construct a partition of $X$. 
Let  $\rho\in \M(G) $ be fixed, for each $h\in H$, put
 	\begin{align}\label{eqn-S1} G_h=\{g\in G:\varphi(g)=h\}\end{align}
 and
 	$$ B_h^{\times}=\bigcup_{k=1}^Kh^{-1}(B_k)\times h^{-1}(B_k).$$
 	For $g\in G$, put
 \begin{align}\label{q-1}
 W_g^\times=\{(x,y)\in B_{\varphi(g)}^\times:d(gx,gy)>2\delta\}.
 \end{align}
 	Then $(x,y)\in W_g^\times$ implies $d(gx,gy)>2\delta$ and $d(\varphi(g)x,\varphi(g)y)<\delta$ (since $\text{diam}(B_k)<\delta,\forall 1\leq k\leq K$). Put
 	\begin{align}\label{15}\begin{split}
 	U_h^\times&:=\{(x,y)\in B_h^{\times}: \rho(\{g\in G_h:(x,y)\in W_g^\times\})\ge \delta \rho(G_h)\}\\
 	&\overset{\eqref{q-1}}=\{(x,y)\in B_h^{\times}: \rho(\{g\in G_h:d(gx,gy)>2\delta\})\ge \delta \rho(G_h)\}.
 	\end{split}
 	\end{align}
We claim that
 	\begin{align}\label{25}(\mu\times\mu)(U_h^\times)\le \frac{\delta^2}{K}.
 	\end{align}
	
	By Cauchy-Schwarz inequality, we have for $h\in H$
 	\begin{align*}
\rho(G_h)\frac{\delta^4}{K}\overset{\eqref{eq-12}}\ge  &	\int_{G_h}\|gf-\varphi(g)f\|_{L^2(\mu\times\mu)}d\rho(g)\\
 		&\ge 		\int_{G_h}\int_{X\times X}|d(gx,gy)-d(\varphi(g)x,\varphi(g)y)|d(\mu\times\mu)(x,y)d\rho(g)\\
 			&\ge 		\int_{G_h}\int_{W_g^\times}|d(gx,gy)-d(\varphi(g)x,\varphi(g)y)|d(\mu\times\mu)(x,y)d\rho(g)\\
 		&\ge		\int_{G_h}\int_{W_g^\times}\delta d(\mu\times\mu)(x,y)d\rho(g)\\
 		&=\delta 	\int_{G_h}(\mu\times\mu)(W_g^\times)d\rho(g).
 	\end{align*}
Then, according to the above inequality, for $h\in H$, one has
 	\begin{align}\label{12}\begin{split}
 	&	\int_{X\times X}\rho(\{g\in G_h:(x,y)\in W_g^\times\})d(\mu\times \mu)(x,y)\\
 	&=\int_{X\times X}\int_{G_h}1_{W_g^\times}(x,y)d\rho(g)d(\mu\times \mu)(x,y)\\
 		&=\int_{G_h}(\mu\times\mu)(W_g^\times)d\rho(g)\\
 		&\le \rho(G_h)\frac{\delta^3}{K}.
 	\end{split}
 	\end{align}	
	By combining \eqref{15} and \eqref{12}, one has the claim proven, i.e.,  \eqref{25} holds.

 Notice that every $U_h^\times$ is a measurable subset of $B_h^{\times}=\bigcup_{k=1}^Kh^{-1}(B_k)\times h^{-1}(B_k)$. By Fubini's Theorem, for each $k\in \{1,\dots, K\}$ with $\mu(h^{-1}B_k)>0$, we can find an $x_{k,h}\in h^{-1}B_k$ and a measurable subset $U_{k,h}$ of $h^{-1}B_k$ such that
 	$$\{x_{k,h}\}\times U_{k,h}=U_h^\times\cap (\{x_{k,h}\}\times h^{-1}B_k)$$
 	and
 	$$\mu(U_{k,h})\le\frac{ \mu\times\mu(h^{-1}(B_k)\times h^{-1}(B_k)\cap U_h^\times)}{\mu(h^{-1}B_k)}. $$
 	Denote
 	$$U_{k,h}^*=h^{-1}B_k\setminus U_{k,h}.$$
	By \eqref{15}, one has \begin{align}\label{eqn-S12}
	U_{k,h}^{*}=\{(x,y)\in B_h^{\times}: \rho(\{g\in G_h:d(gx,gy)>2\delta\})< \delta \rho(G_h)\}. \end{align}
 	Since $\{B_k\}_{k=1}^K$ is a partition of $X$, so is $\{h^{-1}B_k\}_{k=1}^{K}$. Thus for every $h\in G$,
 	\begin{align}\label{q-2}
 	X=\sqcup_{k=1}^K U_{k,h}\bigsqcup\sqcup_{k=1}^K U_{k,h}^*,
 	\end{align}
	which is the desired partition of $X$. 
	
	\smallskip
	
\noindent{\bf Step 2:} In this step, we show that for the given $\mu\in \mathcal{M}(X)$, the fixed $\rho\in \mathcal{M}(G)$, $k\in \{1,\dots, K\}$ and the finite subset $H\subseteq G$ as in the beginning of this proof, one has \begin{align}\label{13}\sum_{h\in H}\rho(G_h)	\mu(\cup_{k=1}^KU_{k,h})\le 2\delta\sum_{h\in H} \rho(G_h),\end{align} where $G_h$ and $U_{k, h}$ are as in \eqref{eqn-S1} and \eqref{q-2}, respectively.

 	For $h\in H$ with $\rho(G_h)>0$ and $x,y\in U_{k,h}^*$, one has
 	\begin{align}\label{16}
 	d_{\rho|_{G_h}}(x,y)\le d_{\rho|_{G_h}}(x,x_{k,h})+d_{\rho|_{G_h}}(x_{k,h},y) 	\le 6\delta.
 	\end{align}
 This comes from the fact that
  \begin{align*}
 	& d_{\rho|_{G_h}}(x,x_{k,h})= \frac{1}{\rho(G_h)}\int_{G_h}d(gx,gx_{k,h})d\rho(g)\\
 	&\le \frac{1}{\rho(G_h)}\int_{\{g\in G_h:d(gx,gx_{k,h})>2\delta\}}d(gx,gx_{k,h})d\rho(g)+\\
 	&\frac{1}{\rho(G_h)}\int_{\{g\in G_h:d(gx,gx_{k,h})\leq2\delta\}}d(gx,gx_{k,h})d\rho(g)\\
 	&\le \frac{1}{\rho(G_h)}\int_{\{g\in G_h:d(gx,gx_{k,h})>2\delta\}}\text{diam}(X)d\rho(g)+2\delta\\
 	&\le \frac{1}{\rho(G_h)}\rho(\{g\in G_h:d(gx,gx_{k,h})>2\delta\})\text{diam}(X)+2\delta\\
 	&\le \delta\text{diam}(X)+2\delta\overset{\eqref{diam1}}=3\delta,
 		\end{align*}
		where the last inequality uses the fact that  $(x,x_{k,h})\in U^{*}_{k,h},\ \forall x\in U_{k,h}^{*}$ and \eqref{eqn-S12}.
		
 	Notice that $\{U_{k,h}\}_{k=1}^K$ are pairwise disjoint. For every $h\in H$, one has
 	\begin{align*}
 		\mu(\cup_{k=1}^KU_{k,h})&=\sum_{k=1}^K\mu(U_{k,h})\\
 		&\le \sum_{1\le k\le K\atop \mu(h^{-1}B_k)<\frac{\delta}{K}}\mu(U_{k,h})+\sum_{1\le k\le K\atop \mu(h^{-1}B_k)\ge\frac{\delta}{K}}\mu(U_{k,h})\\
 		&\le \sum_{1\le k\le K\atop \mu(h^{-1}B_k)<\frac{\delta}{K}}\frac{\delta}{K}+\sum_{1\le k\le K\atop \mu(h^{-1}B_k)\ge\frac{\delta}{K}}\frac{  \mu\times\mu(h^{-1}(B_k)\times h^{-1}(B_k)\cap U_h^\times)}{\mu(h^{-1}B_k)}\\
 		&\le \delta+\frac{ \sum_{1\le k\le K} \mu\times\mu(h^{-1}(B_k)\times h^{-1}(B_k)\cap U_h^\times)}{\delta/K}\\
 			&=\delta+\frac{ \mu\times\mu(U_h^\times)}{\delta/K}\\
 		&\overset{\eqref{25}}\le 2\delta.
 	\end{align*}
 	This implies
 	\begin{align}
 	\int_X\sum_{h\in H} \rho(G_h)	1_{\cup_{k=1}^KU_{k,h}}(x)d\mu(x)=\sum_{h\in H}\rho(G_h)	\mu(\cup_{k=1}^KU_{k,h})\le 2\delta\sum_{h\in H} \rho(G_h),
 	\end{align}
	which finishes the proof of \eqref{13}.
 	\smallskip
	
	\noindent{\bf Step 3:} In this step, for the given $\mu\in \mathcal{M}(X)$, the fixed $\rho \in \mathcal{M}(G)$ and $\delta>0$, we construct a subset $P\subseteq X$ with $\mu(P)>1-\delta$ and can be covered by a finite number of small subsets of $X$. By small subsets, we mean the diameter  of each subset is smaller than $4\sqrt{\delta}+6\delta$ with respect to the metric $d_{\rho}$.   
	
	Set
\begin{align}\label{18}\begin{split}
P&=\{x\in X: \sum_{h\in H} \rho(G_h)	1_{\cup_{k=1}^KU_{k,h}}(x)<2\sqrt{\delta}\sum_{h\in H} \rho(G_h)\}\\
&\overset{\eqref{q-2}}=\{x\in X: \sum_{h\in H} \rho(G_h)	1_{\cup_{k=1}^KU^*_{k,h}}(x)> (1- 2\sqrt{\delta})\sum_{h\in H} \rho(G_h)\}.
\end{split}
\end{align}
Clearly, for the complement of $P$, i.e., $$P^{c}:=\{x\in X: \sum_{h\in H}\rho(G_{h}I_{\cup_{k=1}^{K}U_{k,h}})(x)\geq 2\sqrt{\delta}\sum_{h\in H}\rho(G_{h})\}.$$ 
By  \eqref{13}, one has $ \mu(P^c)\leq \sqrt{\delta}$.
 	Then, we obtain
 \begin{align}\label{19}\mu(P)\ge 1-\sqrt{\delta}.\end{align}
 
 	Now, we demonstrate that $P\subseteq X$ with $\mu(P)>1-\delta$ can be covered by a finite number of small subsets. For any sequence of length $|H|$ with each element chosen from $\{1,\dots, K\}$, i.e., $(k_h)_{h\in H}\in \{1,2,\dots,K\}^H$, we put
 \begin{align}\label{20}
 P((k_h)_{h\in H})=\{x\in X:\sum_{h\in H} \rho(G_h)	1_{U_{k_h,h}^*}(x)\ge (1- 2\sqrt{\delta})\sum_{h\in H} \rho(G_h)\}.
 \end{align}
 By \eqref{q-2},  $	1_{\cup_{k=1}^KU^*_{k,h}}(x)=\sum_{k=1}^K	1_{U^*_{k,h}}(x)$. Combining this with \eqref{18} and \eqref{20}, one has
 $$P=\bigcup_{(k_h)_{h\in H}\in \{1,2,\dots,K\}^H}P((k_h)_{h\in H}).$$
 Thus
 \begin{align}\label{14}
 \mu\left(\bigcup_{(k_h)_{h\in H}\in \{1,2,\dots,K\}^H}P((k_h)_{h\in H})\right)=\mu(P)\overset{\eqref{19}}\ge 1-\sqrt{\delta}.
 \end{align}
For $(k_h)_{h\in H}\in \{1,2,\dots,K\}^H$ and $x,y\in P((k_h)_{h\in H})$, we deduce that
 	\begin{align*} d_\rho(x,y)&=\sum_{h\in H, \rho(G_h)>0} \rho(G_h) d_{\rho|_{G_h}}(x,y)\\
 		&=\sum_{h\in H,\rho(G_h)>0\atop x\notin U_{{k_h},h}^*\text{or }y\notin U_{{k_h},h}^*} \rho(G_h) d_{\rho|_{G_h}}(x,y)+\sum_{h\in H,\rho(G_h)>0\atop x,y\in U_{{k_h},h}^* } \rho(G_h) d_{\rho|_{G_h}}(x,y)\\
 		&\overset{\eqref{16}}\le \sum_{h\in H,\rho(G_h)>0\atop x\notin U_{{k_h},h}^*\text{or }y\notin U_{{k_h},h}^*} \rho(G_h) \text{diam}(X)+\sum_{h\in H,\rho(G_h)>0\atop x,y\in U_{{k_h},h}^* } \rho(G_h)(6\delta)\\
 		&\le \sum_{h\in H\atop x\notin U_{{k_h},h}^*} \rho(G_h)\text{diam}(X)+\sum_{h\in H\atop y\notin U_{{k_h},h}^*} \rho(G_h)\text{diam}(X)+6\delta\\
 		&\overset{\eqref{20}}\le 4\sqrt{\delta}+6\delta.
 	\end{align*}
 Combining the above inequality with \eqref{14}, one has
 	$$S_\rho(X,\mu,G,6\delta+4\sqrt{\delta})\le K^{|H|}.$$
This implies $\mu$ has bounded measure-max-mean-complexity, which completes the proof of (1)$\Rightarrow$(2). 	

 \smallskip	
 		
 \noindent{\bf (3)$\Rightarrow$(1):} We assume that the  max-mean-complexity of  $\mu$ is bounded.  For a contradiction, assume that $\mu\in\mathcal{M}(X)$ does not have discrete spectrum. We can find a continuous function $f$ on $X$, an infinite subset $\mathcal{S}$ of $G$ and a constant $\epsilon>0$ such that for distinct $g,h\in \mathcal{S}$ one has
 	\begin{align}\label{equition_1}\|hf-gf\|_{L^2(\mu)}\ge \epsilon.\end{align}
 	Now we fix $\epsilon',\delta>0$ sufficiently small such that the followings hold.
 	\begin{itemize}
 		\item[(C1).]If $x,y\in X,d(x,y)\le \delta$, then  $|f(x)-f(y)|\le \epsilon'$.
 		\item[(C2).] $4\|f\|_\infty(6(\epsilon'+2\delta\|f\|_\infty)+2\delta^2\|f\|_\infty)
 		<\epsilon^2.$
 	\end{itemize}

Recall that $\rho_E=\frac{1}{|E|}\sum_{g\in E}\delta_g$ for $E\in F(G)$. Since the  max-mean-complexity of  $\mu$ is bounded, there exists $K\in\N$ such that
 	$$S_{\rho_E}(X,\mu,G,\delta^2)\le K$$ for any finite subset $E$ of $\mathcal{S}$.
 	Let $\{C_j: j\in\{1,2,\dots,J\}\}$ be a finite Borel measurable partition of $X^K$ with $\text{diam}(C_j)<\delta,\forall j\in\{1,2,\dots,J\}$, that is, if $(y_1,y_2,\dots,y_K)$, $(y_1',y_2',\dots,y_K')\in C_j$ then
 	$$d(y_k,y_k')<\delta,\ \forall k\in\{1,2,\dots,K\}.$$
 	Fix a finite subset $E$ of $\mathcal{S}$ such that
 	\begin{align}\label{1}|E|>2J.
 	\end{align}  By the definition of $S_{\rho_E}(X,\mu,G,\delta^2)$ we can find  $x_1,x_2,\dots,x_K\in X$ such that
 		\begin{align}\label{q-4}\mu(\bigcup_{k=1}^K B_{\rho_E}(x_i,\delta^2))>1-\delta^2.\end{align}	 For each $j\in\{1,2,\dots,J\}$, put
 	$$E_j=\{g\in E:(gx_1,gx_2,\dots,gx_K)\in C_j\}.$$
 	By \eqref{1}, one has
 	\begin{align}\label{2}\begin{split}
 			|E|&=2|E|-|E|=2\left(\sum_{j\in\{1,2,\dots,J\}\atop |E_j|=1}1+\sum_{j\in \{1,2,\dots,J\}\atop |E_j|\ge 2}|E_j|\right)-|E|\\
 			&\le 2\left(J+\sum_{j\in \{1,2,\dots,J\}\atop |E_j|\ge 2}|E_j|\right)-|E|\overset{ \eqref{1} }\le |E|+2\sum_{j\in \{1,2,\dots,J\}\atop |E_j|\ge 2}|E_j|-|E|\\
 			&=2\sum_{j\in\{1,2,\dots,J\}\atop |E_j|\ge 2}|E_j|.
 		\end{split}
 	\end{align}
 	For any $k\in\{1,2,\dots,K\}$ and $x\in B_{\rho_E}(x_k,\delta^2)$, i.e., $\frac{1}{|E|}\sum_{g\in E} d(gx, gx_k)<\delta^2$, we obtain
 	$$\frac{|\{g\in E:d(gx,gx_k)\geq \delta\}|}{|E|}<\delta.$$
 	Combining with (C1), one has
 	\begin{align}\label{eqn-3}\frac{1}{|E|}\sum_{g\in E}|f(gx)-f(gx_k)|
	\le \epsilon'+2\delta\|f\|_\infty.\end{align}
 	Therefore, for each $k\in \{1,\dots, K\}$ and $x\in B_{\rho_E}(x_k, \delta^2)$, one has
 	\begin{align}\label{q-3}\begin{split}
 	&\sum_{j\in\{1,2,\dots,J\}\atop
 		|E_j|\ge 1}\frac{1}{|E_j|}\sum_{g,h\in E_j}|f(hx)-f(gx)|\\
 	&\le \sum_{j\in\{1,2,\dots,J\}\atop
 		|E_j|\ge 1}\frac{1}{|E_j|}\sum_{g,h\in E_j}\left(|f(hx)-f(hx_k)|+|f(hx_k)-f(gx_k)|+|f(gx_k)-f(gx)|\right)\\
 	&=  2\sum_{j\in\{1,2,\dots,J\}\atop
 		|E_j|\ge 1}\sum_{g\in E_j}|f(gx)-f(gx_k)|+\sum_{j\in\{1,2,\dots,J\}\atop
 		|E_j|\ge 1}\frac{1}{|E_j|}\sum_{g,h\in E_j}|f(hx_k)-f(gx_k)|\\
 	&\le 2\sum_{g\in E}|f(gx)-f(gx_k)|+\sum_{j\in\{1,2,\dots,J\}\atop
 		|E_j|\ge 1}\frac{1}{|E_j|}\sum_{g,h\in E_j}\epsilon'\\
 	&\overset{\eqref{eqn-3}}{\le} 3|E|(\epsilon'+2\delta\|f\|_\infty).
 	\end{split}
 	\end{align}
 	Put $W=\bigcup_{k=1}^K B_{\rho_E}(x_k,\delta^2)$. Then
 	\begin{align*} &\sum_{j\in\{1,2,\dots,J\}\atop
 			|E_j|\ge 2}|E_j|\frac{1}{|E_j|^2}\sum_{g,h\in E_j}\int_{X}1_W(x)|f(hx)-f(gx)|d\mu(x)\\
 		&\le \int_{X}\sum_{j\in\{1,2,\dots,J\}\atop
 			|E_j|\ge 1}\frac{1_W(x)}{|E_j|}\sum_{g,h\in E_j}|f(hx)-f(gx)|d\mu(x)\\
 		&\overset{\eqref{q-3}}\le 3| E|(\epsilon'+2\delta\|f\|_\infty)\\
 		&\overset{\eqref{2}}\le 6\sum_{j\in\{1,2,\dots,J\}\atop
 			|E_j|\ge 2} |E_j|(\epsilon'+2\delta\|f\|_\infty).
 	\end{align*}
 	We can find $s\in\{1,2,\dots,J\}$ such that
 	$$|E_s|\ge 2$$
 	and
 	$$\frac{1}{|E_s|^2}\sum_{g,h\in E_s}\int_{X}1_W(x)|f(hx)-f(gx)|d\mu(x)\le6(\epsilon'+2\delta\|f\|_\infty).$$
 	This implies
 	\begin{align*}
 		\frac{1}{|E_s|^2}\sum_{g,h\in E_s}\int_{X}|f(hx)-f(gx)|d\mu(x)&\le6(\epsilon'+2\delta\|f\|_\infty)+2\mu(X\setminus W)\|f\|_\infty\\
 		&\overset{\eqref{q-4}}\le 6(\epsilon'+2\delta\|f\|_\infty)+2\delta^2\|f\|_\infty.
 	\end{align*}
 Since $|E_s|\ge 2$,	we can find distinct $g,h\in E_s\subseteq \mathcal{S}$ such that
 $$\int_{X}|f(hx)-f(gx)|d\mu(x)\le 2(6(\epsilon'+2\delta\|f\|_\infty)+2\delta^2\|f\|_\infty).$$
 Thus
 	\begin{align*}\|hf-gf\|_{L^2(\mu)}^2&\le 2\|f\|_\infty\int_{X}|f(hx)-f(gx)|d\mu(x)\\
 		&\le 4\|f\|_\infty(6(\epsilon'+2\delta\|f\|_\infty)+2\delta^2\|f\|_\infty)\\
 		&\overset{(C2)}<\epsilon^2.
 	\end{align*}
 	However, this yields a contradiction to \eqref{equition_1}. We end the proof of Theorem \ref{thm-2}.
 \end{proof}
\section{Proof of Theorem \ref{thm-1}}\label{s-t-2}
 In this section, we are to show that for locally compact amenable group action $(X, G)$, the discrete spectrum, the bounded max-mean complexity along any F\o lner sequence of $G$ and the bounded measure max-mean-complexity of any $G$-invariant Borel probability measure are equivalent.
 \begin{proof}[Proof of Theorem \ref{thm-1}](3)$\Leftrightarrow$(4) is from Theorem \ref{thm-2}. (3)$\Rightarrow$(2)$\Rightarrow$(1)  obviously holds by the definitions of  the bounded measure-max-mean-complexity and the bounded mean-complexity along a F\o lner sequence of $G$. The only left is to prove (1)$\Rightarrow$(4). 

Let $G$ be a locally compact amenable group, $m_G$ be a (left) Haar measure on $G$ and $\mu$ be a $G$-invariant Borel probability measure on $X$. Without loss of generality, we assume that
 	 \begin{align}\label{eq-1}\text{diam}(X)=1.
 	\end{align}
 	We assume that $\mu$ has bounded mean-complexity along the F{\o}lner sequence $\{F_n\}_{n\in\N}$ of $G$. For a given $\epsilon>0$, by the definition of bounded mean-complexity along the F{\o}lner sequence $\{F_n\}_{n\in\N}$, there exists $K=K(\epsilon)\in \mathbb{N}$ such that
\begin{align}\label{e-0} S_{m_G|_{F_n}}(X,\mu,G,\epsilon)<K,\forall n\in\N.\end{align}
 	To show $\mu\in\mathcal{M}(X)$ has discrete spectrum, by Theorem \ref{thm-2} it is sufficient to prove that $\mu$ has bouned max-mean-complexity, i.e.,
 	 $$ S_{\rho_E}(X,\mu,T,2\sqrt{\epsilon})\le K$$
 	 for every finite subset $E$ of $G$ where  $\rho_E=\frac{1}{|E|}\sum_{g\in E}\delta_g$.

 Given a finite subset $E$ of $G$. Fix $n\in\N$ sufficiently large such that
 	\begin{align}\label{q-6}\frac{m_G((EF_n)\triangle F_n)}{m_G(F_n)}\le \frac{\epsilon}{|E|}.\end{align}
 	Put
 \begin{align}\label{eq-0}H=\{h\in F_n:Eh\subseteq F_n\}.
 \end{align}
 	If $h\in F_n\setminus H$ then $Eh\setminus F_n\neq\emptyset$ which implies $h\in E^{-1}((EF_n)\setminus F_n)$.
 	Thus
 	$$m_G(H)\ge m_{G}(F_n)-m_G(E^{-1}((EF_n)\setminus F_n))\overset{\eqref{q-6}}\ge (1-\epsilon)m_G(F_n).$$ 	Since $m_G$ is left invariant,
 	\begin{align*}
 m_G(gH)= m_G(H)\ge  (1-\epsilon)m_G(F_n),\text{ for all }g\in E.
 	\end{align*}
 	Since $gH\subseteq F_n$ for all $g\in E$, one has
 	\begin{align}\label{eq-123}
 m_G|_{F_n}(gH)\ge 1-\epsilon, \text{ for all }g\in E.
 	\end{align}
 	By \eqref{e-0}, there exist $x_1,x_2,\dots,x_K$ of $X$ such that
 	\begin{align}\label{q-7}\mu(\bigcup_{k=1}^K  B_{m_G|_{F_n}}(x_k,\epsilon))>1-\epsilon.\end{align}
 	For each $k\in\{1,2,\dots,K\}$ and $x\in  B_{m_G|_{F_n}}(x_k,\epsilon)$, one has
 \begin{align}\label{22}\begin{split}\int_H\sum_{g\in Eh}d(gx,gx_k)dm_G|_{F_n}&(h)=\sum_{g\in E}\int_Hd(ghx,ghx_k)dm_G|_{F_n}(h)\\
 &\overset{H\subseteq F_n}= \frac{1}{m_G(F_n)}\sum_{g\in E}\int_Hd(ghx,ghx_k)dm_G(h)\\
 &=\frac{1}{m_G(F_n)}\sum_{g\in E}\int_{gH}d(hx,hx_k)dm_G(h)\\
 &\overset{EH\subseteq F_n}=\sum_{g\in E}\int_{gH}d(hx,hx_k)dm_G|_{F_n}(h)\\
  &\overset{\eqref{eq-123}}\le \sum_{g\in E}\left(\int_{F_n}d(hx,hx_k)dm_G|_{F_n}(h)+\epsilon\text{diam}(X)\right)\\
 &\le 2|E|\epsilon.
 \end{split}
 \end{align}
 	For each $k\in\{1,2,\dots,K\}$ and $x\in   B_{m_G|_{F_n}}(x_k,\epsilon)$, set
 \begin{align}\label{123}
 \mathcal{C}_k(x)=\{h\in H:\sum_{g\in Eh}d(gx,gx_k)\le |E|2\sqrt{\epsilon}\}.
 \end{align}
 	Then by \eqref{22} for each $k\in\{1,2,\dots,K\}$ and $x\in   B_{m_G|_{F_n}}(x_k,\epsilon)$,
 	\begin{align}\label{eq-02}m_G|_{F_n}(\mathcal{C}_k(x))\ge 1-\sqrt{\epsilon}.
 	\end{align}
 	Put $B_1= B_{m_G|_{F_n}}(x_1,\epsilon)$ and $B_k=B_{m_G|_{F_n}}(x_1,\epsilon)\setminus \cup_{l=1}^{k-1}B_{m_G|_{F_n}}(x_l,\epsilon)$ for $k\in\{2,\dots,K\}$.

 Thus
 	\begin{align}\label{eq-03}\mu(\cup_{k=1}^KB_k)=\mu(\cup_{k=1}^KB_{m_G|_{F_n}}(x_k,\epsilon))\overset{\eqref{q-7}}\ge 1-\epsilon.
 	\end{align} 	For each $k\in\{1,2,\dots,K\}$ and $h\in H$, put
 	$$Y_{k,h}=\{x\in  B_k:\sum_{g\in Eh}d(gx,gx_k)\le |E|2\sqrt{\epsilon}\}.$$
 	Thus, we get
 	\begin{align}\label{222}
  Y_{k,h}\subseteq B_{\rho_{Eh}}(x_k,2\sqrt{\epsilon}).
 	\end{align}
 	By \eqref{123}, we know
 \begin{align}\label{eq-01}x\in Y_{k,h}\text{ if and only if }x\in B_k\text{ and }h\in\mathcal{C}_k(x).\end{align}
 	Thus by Fubini's Theorem
 \begin{align*}
 \int_H\sum_{k=1}^K\mu(Y_{k,h})dm_G|_{F_n}(h)&\overset{\eqref{222}}{=} \int_H\sum_{k=1}^K\int_{ B_k} 1_{Y_{k,h}}(x)d\mu(x)dm_G|_{F_n}(h)\\
 &= \sum_{k=1}^K\int_{ B_k}\int_H 1_{Y_{k,h}}(x)dm_G|_{F_n}(h)d\mu(x)\\
  &\overset{\eqref{eq-01}}= \sum_{k=1}^K\int_{ B_k}m_G|_{F_n}( \mathcal{C}_k(x))d\mu(x)\\
  &\overset{\eqref{eq-02}}\ge  \sum_{k=1}^K\int_{ B_k}(1-\sqrt{\epsilon})d\mu(x)\\
 &\overset{\eqref{eq-03}}\ge (1-\epsilon)(1-\sqrt{\epsilon}).
 \end{align*}
 	There is $h\in H$ such that
 	$$\sum_{k=1}^K\mu(Y_{k,h})\ge (1-\epsilon)(1-2\sqrt{\epsilon})> 1-2\sqrt{\epsilon}.$$
 Since $\{Y_{k,h}\}_{k=1}^K$ are pairwise disjoint, one has
 $$\mu(\cup_{k=1}^K B_{\rho_{Eh}}(x_k,2\sqrt{\epsilon}))\overset{\eqref{222}}\ge \sum_{k=1}^K\mu(Y_{k,h})> 1-2\sqrt{\epsilon}.$$ Therefore, $$ S_{\rho_E}(X,\mu,G,2\sqrt{\epsilon})=  S_{\rho_{Eh}}(X,\mu,G,2\sqrt{\epsilon})\le K.$$
 	By the arbitrariness of $E$,  $\mu$ has bounded max-mean-complexity. By Theorem \ref{thm-2}, $\mu$ has discrete spectrum. This ends the proof of Theorem \ref{thm-1}.
 \end{proof}
 
 \section{Proof of Theorem \ref{thm-3}}\label{s-t-3}
In this section, the proof of Theorem \ref{thm-3} will be given. Specifically, for a locally compact amenable group action $(X, G)$, the discrete spectrum and mean equicontinuity along a tempered F\o lner sequence (equicontinuity in the mean along a tempered F\o lner sequence) of any $G$-invariant Borel probability measure are equivalent. 

Firstly, we need the following lemma, which gives an equivalent definition of mean equicontinuous along  a F{\o}lner sequence and equicontinuity in the mean along a F\o lner sequence.
\begin{lem}\label{lem-2}Let $G$ be a locally compact amenable group and $(X,\mathcal{B}_X,\mu,G)$ be a measure-preserving system with a F{\o}lner sequence $\{F_n\}_{n\in \mathbb{N}}.$   Then the followings hold.
	\begin{itemize}
		\item[(1)] $\mu$ is mean equicontinuous along $\{F_n\}_{n\in  \mathbb{N}}$ if and only if for any $\epsilon>0$ there exists  a $\delta>0$ and a measurable set $Q\subseteq X$ with $\mu(Q)>1-\epsilon$ such that 
		$$x,y\in Q, d(x,y)<\delta \Rightarrow \limsup_{n\to\infty}d_{m_G|_{F_n}}(x,y)\le 2\epsilon.$$
		\item[(2)] $\mu$ is equicontinuous in the
		mean  along $\{F_n\}_{n\in  \mathbb{N}}$ if and only if for any $\epsilon>0$ there exists a $\delta>0$ and a measurable set $Q\subseteq X$ with $\mu(Q)>1-\epsilon$ such that 
		$$x,y\in Q, d(x,y)<\delta \Rightarrow d_{m_G|_{F_n}}(x,y)\le 2\epsilon,\ \forall n\in \mathbb{N}.$$
	\end{itemize}
\end{lem}
\begin{proof}(1)
		 ``$\Rightarrow$''  is trivial. We prove the opposite side. By assumption, for $j\in\N$, we can find a $\delta_j>0$ and a measurable set $Q_j\subseteq X$ with $\mu(Q_j)>1-1/j^2$ such that 
		$$x,y\in Q_j, d(x,y)<\delta_j \Rightarrow \limsup_{n\to\infty}d_{m_G|_{F_n}}(x,y)<2/j^2.$$
		For given $\tau>0$, we can find a $J\in\N$ such that $\sum_{j\ge J}1/j^2<\tau$. Put $Q=\cap_{j\ge J} Q_j$. Then
		$$\mu(Q)\ge 1-\tau.$$
		And for every $\epsilon>0$, we can find  a $j_\epsilon\ge J$ such that $2/j_\epsilon^2<\epsilon$. Let $\delta=\delta_{j_\epsilon}$. Since $Q\subseteq Q_{j_\epsilon}$, one has  
		$$x,y\in Q, d(x,y)<\delta \Rightarrow \limsup_{n\to\infty}d_{m_G|_{F_n}}(x,y)<2/j_\epsilon^2<\epsilon.$$
		Thus,  $\mu$ is mean equicontinuous along $\{F_n\}_{n\in  \mathbb{N}}$.\
		
		(2)	 ``$\Rightarrow$''  is trivial. We prove the opposite side. By assumption, for $j\in\N$, we can find a $\delta_j>0$ and a measurable set  $Q_j\subseteq X$ with $\mu(Q_j)>1-1/j^2$ such that 
		$$x,y\in Q_j, d(x,y)<\delta_j \Rightarrow d_{m_G|_{F_n}}(x,y)<2/j^2, \forall n\in \mathbb{N}.$$
		For given $\tau>0$, we can find $J\in\N$ such that $\sum_{j\ge J}1/j^2<\tau$. Put $Q=\cap_{j\ge J} Q_j$. Then
		$$\mu(Q)\ge 1-\tau.$$
	And for every $\epsilon>0$, we can find  a $j_\epsilon\ge J$ such that $2/j_\epsilon^2<\epsilon$. Let $\delta=\delta_{j_\epsilon}$. Since $Q\subseteq Q_{j_\epsilon}$, one has  
	$$x,y\in Q,\ d(x,y)<\delta \Rightarrow d_{m_G|_{F_n}}(x,y)<2/j_\epsilon^2<\epsilon, \forall n\in \mathbb{N}.$$
		Thus,  $\mu$ is equicontinuous in the
		mean  along $\{F_n\}_{n\in  \mathbb{N}}$.
\end{proof}
%
Now, we are ready to prove Theorem \ref{thm-3}.
\begin{proof}[Proof of Theorem \ref{thm-3}]
There are three parts of the proof, i.e.,
\begin{itemize}
\item[(1)$\Rightarrow$(2):]  $\mu\in\mathcal{M}(X, G)$ has discrete spectrum implies it is mean equicontinuous along a tempered F\o lner sequence $\{F_n\}_{n\in\mathbb{N}}$;
\item[(2)$\Rightarrow$(3):] If $\mu\in \mathcal{M}(X, G)$ is mean equicontinuous along a tempered F\o lner sequence $\{F_n\}_{n\in\mathbb{N}}$, then it is equicontinuous in the mean along the same F\o lner sequence; 
\item[(3)$\Rightarrow$(1):]  $\mu\in \mathcal{M}(X,G)$ is equicontinuous in the mean along a F\o lner sequecne $\{F_n\}_{n\in \mathbb{N}}$ shows that it has discrete spectrum.
\end{itemize}
\smallskip

\noindent{\bf (1)$\Rightarrow$(2):} By Theorem \ref{thm-1}, we know that $\mu$ has bounded max-mean-complexity along the tempered F{\o}lner sequence $\{F_n\}_{n\in \mathbb{N}}$. Given $0<\epsilon<\frac{1}{2}.$ Then there exists $K\in\N$, and for each $n\in \mathbb{N}$, there is a family of measurable subsets $\{C_1^n,C_2^n,\dots,C_K^n\}$ of $X$ such that for $n\in\N$
		$$\mu(\cup_{k=1}^K C_k^n)> 1-\epsilon/2,$$
		and
			\begin{align}\label{e-6} d_{m_G|_{F_n}}(x,y)<\epsilon\text{ for }x,y\in C_k^n, k\in\{1,2,\dots,K\}.\end{align}
			 Additionally, we can assume that each $C_k^n$ has measure larger than $\frac{\epsilon}{2K}$ (If the measure of some $C_k^n$ is smaller than $\frac{\epsilon}{2K}$, we can change it to the bigger one). Thus it follows
		\begin{align}\label{e-1}(\mu\times\cdots\times\mu)(\cap_{N\in\N}\cup_{n\ge N}C_1^n\times\cdots\times C_K^n)>(\frac{\epsilon}{2K})^K.
		\end{align}
		Let $f(x,y)=d(x,y)$ for $x,y\in X$. Let
			\begin{align}\label{e-2} X_*^\times=\{ (x,y)\in X\times X: \lim_{n\to\infty}\frac{1}{m_G(F_n)}\int_{F_n} f(gx,gy)dm_G(g)\text{ exists}\}.\end{align}
		Then by the pointwise ergodic theorem for amenable group actions \cite{Li}, $\mu\times\mu(X_*^\times)=1$. For $x\in X$, put
		$$R_x=\{y\in X:(x,y)\in X_*^\times\}.$$ We can choose $(x_1,\dots, x_K)\in \cap_{N\in\N}\cup_{n\ge N}C_1^n\times\cdots\times C_K^n$ such that
		\begin{align}\label{e-3}\mu(R_{x_k})=1\text{ for }k=1,2,\dots,K.\end{align}
	There exists $n_1<n_2<\cdots$ such that
	$$x_k\in C_{k}^{n_i}\text{ for }k=1,2,\dots,K\text{ and }i\in \mathbb{N}.$$
	Let
		\begin{align}\label{e-4}C_*= \left(\cap_{I\in \N}\cup_{i\ge I}\cup_{k=1}^KC_k^{n_i}\right)\bigcap \left(\cap_{k=1}^KR_{x_k}\right).\end{align}
	Then by \eqref{e-1} and \eqref{e-3}, $$\mu(C_*)> 1-\epsilon.$$
	
	We claim that for  every $x\in C_*$, there exists $k_*$ depending on $x$ such that 	\begin{align}\label{e-7}\limsup_{n\to\infty}d_{m_G|_{F_n}}(x,x_{k_*})\le \epsilon.\end{align}
	The proof is postponed to the end of this step. 

	Let
	$$D_k=\{x\in X: \limsup_{n\to\infty}d_{m_G|_{F_n}}(x,x_{k})\le \epsilon\}.$$ By \eqref{e-7},
	$$\mu(\cup_{k=1}^K D_k)\ge \mu(C_*)> 1-\epsilon.$$
	One can find pairwise disjoint compact subsets $E_1\subseteq D_1,E_2\subseteq D_2,\dots, E_K\subseteq D_K$ such that
	\begin{align}\label{e-8}\mu(\cup_{k=1}^K E_k)> 1-\epsilon.
	\end{align}
	Without loss of generality, assume that each $E_k$ is not empty.  Let $\delta>0$ be  the minimum mutual distance between $E_1,\dots,E_K$. Whenever $x,y\in \cup_{k=1}^KE_k, d(x,y)<\delta$, one can find $k\in\{1,2,\dots,K\}$ such that $x,y\in E_k\subseteq D_k$. Then
	$$\limsup_{n\to\infty}d_{F_n}(x,y)\le \limsup_{n\to\infty}d_{F_n}(x,x_k)+\limsup_{n\to\infty}d_{F_n}(y,x_k)\le 2\epsilon.$$
	It follows that
	$$x,y\in \cup_{k=1}^KE_k, d(x,y)<\delta\text{ implies }\limsup_{n\to\infty}d_{F_n}(x,y)\le 2\epsilon.$$
 Combining this with \eqref{e-8} and Lemma \ref{lem-2} (1), we have shown (1)$\Rightarrow$(2).	
 	
Now, we prove the claim.
Let $x\in C_*$. There exists $i_1<i_2<\cdots$ depending on $x$ such that for $j\in\N$
	$$x\in \cup_{k=1}^KC_{k}^{n_{i_j}}.$$ By \eqref{e-4}, there exists $k_j\in \{1,2,\dots,K\}$ such that
	$$x\in C_{k_j}^{n_{i_j}}.$$
	Passing to  a subsequence, we can assume that $k_j=k_*,\forall j\in\mathbb{N}$. Then, by \eqref{e-7}, for $x\in C_*$,
	\begin{align}\label{e-89}\limsup_{j\to\infty}\frac{1}{m_G(F_{n_{i_j}})}\int_{F_{n_{i_j}}} f(gx,gx_{k_*})dm_G(g)\le \epsilon.
	\end{align}
	Since $(x_{k_*},x)\in X_*^\times$, by \eqref{e-2}, $\lim_{n\to\infty}\frac{1}{m_G(F_n)}\int_{F_n} f(gx,gx_{k_*})dm_G(g)$ exists. By \eqref{e-89}, for every $x\in C_*$ $$\lim_{n\to\infty}\frac{1}{m_G(F_n)}\int_{F_n} f(gx,gx_{k_*})dm_G(g)\le \epsilon.$$ This implies
	$$\limsup_{n\to\infty}d_{m_G|_{F_n}}(x,x_{k_*})\le \epsilon.$$ 
	We finish the proof of the claim, and the step of (1)$\Rightarrow$(2).
 	\smallskip

 		\noindent{\bf (2)$\Rightarrow$(3):} Assume that $\mu$ is mean equicontinuous along $\{F_n\}_{n\in \mathbb{N}}$.  By Lemma \ref{lem-2} (1), for given $\epsilon\in (0,1)$, there exists a $\delta>0$ and a measurable $Q\subseteq X$ with $\mu(Q)>1-\epsilon/2$ such that 
 		\begin{align}\label{eq-11-1}
 		x,y\in Q, d(x,y)<\delta \Rightarrow \limsup_{n\to\infty}d_{m_G|_{F_n}}(x,y)\le \epsilon/2.
 		\end{align} There exists $K\in\N$ and  a measurable partition $C_1,\dots,C_{K}$ of $Q$ with diameter smaller than $\delta$ and $\mu(C_k)>0$ for $1\le k\le K$. Then 
 		\begin{align}\label{eq-11-15-2}
 		\gamma:=(\mu\times\mu)(\cup_{k=1}^K (C_k\times C_k))>0,
 		\end{align}
 		and 
 		\begin{align}\label{eq-11-15-3}\eta:=\min_{1\le k\le K}\mu(C_k)>0.\end{align}
 		
 		For $N\in\N$, define
 		\begin{align}\label{eq-11-15-7-1}Q_{N}^\times=\{(x,y)\in\cup_{k=1}^K (C_k\times C_k):d_{m_G|_{F_n}}(x,y)<\epsilon, \forall n\ge N\}.\end{align}
 		Then by \eqref{eq-11-1},
 		$$\lim_{N\to\infty}(\mu\times\mu)(Q_{N}^\times)=(\mu\times\mu)(\cup_{k=1}^K (C_k\times C_k))=\gamma.$$
 		We can find $N\in\N$ such that
 		 		 \begin{align}\label{eq-11-15-5}(\mu\times\mu)(Q_{N}^\times)>(1-\epsilon\eta^2/4)\gamma.\end{align}
 		
 		Let  $I$ be the set of $k\in\{1,\dots,K\}$ such that
 		 \begin{align}\label{eq-11-15-8}
 		(\mu\times\mu)((C_k\times C_k)\cap Q_{N}^\times )\ge (1-\epsilon/2) \cdot(\mu\times\mu)(C_k\times C_k),
 		\end{align}
		and for $k\in \{1,\dots, K\}\setminus I$, one has
		\begin{align}\label{eq-11-15-7}
	(\mu\times\mu)((C_k\times C_k)\cap Q_{N}^\times )<  (1-\epsilon/2) \cdot (\mu\times\mu)(C_k\times C_k).
	\end{align}
 	We claim that \begin{align}\label{claim4-1}I=\{1,2,\dots,K\}.\end{align}

By assuming the claim, \eqref{eq-11-15-8} holds for all $k\in\{1,\dots,K\}$. 
 		For $k\in\{1,\dots,K\}$, by Fubini Theorem and \eqref{eq-11-15-8}, there exists $x_k\in C_k$ and a measurable set $D_k\subseteq C_k$ such that 
 		 \begin{align}\label{eq-11-15-9}\{x_k\}\times D_k\subseteq  Q_{N}^\times,\end{align} and 
 		$$\mu(D_k)\ge \frac{(\mu\times\mu)((C_k\times C_k)\cap Q_{N}^\times )}{(\mu\times\mu)(C_k\times C_k)}\mu(C_k)\ge (1-\epsilon/2)\mu(C_k).$$
 	Thus,
 	$$\mu(\cup_{k=1}^KD_k)\ge (1-\epsilon/2)\mu(Q)>1-\epsilon.$$	
 	One can find pairwise disjoint compact subsets $E_1\subseteq D_1,\dots, E_K\subseteq D_K$ such that
 \begin{align}\label{eq-11-15-6}\mu(\cup_{k=1}^K E_k)> 1-\epsilon.\end{align}
 Without loss of generality, assume that each $E_k$ is not empty.  Let $\delta'>0$ be the minimum mutual distance between $E_1,\dots,E_K$. Whenever $x,y\in \cup_{k=1}^KE_k$ with $d(x,y)<\delta'$, one can find $k\in\{1,2,\dots,K\}$ such that $x,y\in E_k\subseteq D_k$.  By \eqref{eq-11-15-7-1} and \eqref{eq-11-15-9},
 $$d_{F_n}(x,y)\le d_{F_n}(x,x_k)+d_{F_n}(x_k,y)\le 2\epsilon \text{ for } n\ge N.$$
 Thus,
 $$x,y\in \cup_{k=1}^KE_k, d(x,y)<\delta'\text{ implies } d_{F_n}(x,y)\le 2\epsilon \text{ for } n\ge N.$$
Let $\delta''>0$ such that 
$$d(x,y)<\delta''\text{ implies } d_{F_n}(x,y)\le 2\epsilon\text{ for } 1\le n\le N-1.$$
Such $\delta''$ exists since $\cup_{n=1}^{N-1}F_n\subseteq G$ is bounded.
Put $\delta'''=\min\{ \delta',\delta''\}$. Then 
$$x,y\in \cup_{k=1}^KE_k, d(x,y)<\delta'''\text{ implies } d_{F_n}(x,y)\le 2\epsilon, \forall n\in\N.$$
 Combining this with \eqref{eq-11-15-6} and Lemma \ref{lem-2} (2), we prove (2)$\Rightarrow$(3).	

Now, we prove \eqref{claim4-1}. By combining \eqref{eq-11-15-5} and \eqref{eq-11-15-7}, one has
	\begin{align*}
 		(1-\epsilon\eta^2/4)\gamma&\le (\mu\times\mu)(Q_{N}^\times)\\
 		&\le \sum_{k\in I}(\mu\times\mu)((C_k\times C_k)\cap Q^{\times}_N)+ \sum_{k\in\{1,\dots,K\}\setminus I}(\mu\times\mu)((C_k\times C_k)\cap Q^{\times}_N)\\
		&\le \sum_{k\in I}(\mu\times\mu)((C_k\times C_k)\cap Q^{\times}_N)+ \sum_{k\in\{1,\dots,K\}\setminus I}(1-\frac{\epsilon}{2})(\mu\times\mu)(C_k\times C_k)\\
		 &=\sum_{k\in\{1,\dots,K\}}(\mu\times\mu)(C_k\times C_k)-\frac{\epsilon}{2}\sum_{k\in\{1,\dots,K\}\setminus I}(\mu\times\mu)(C_k\times C_k)\\
		 &=\gamma-\frac{\epsilon}{2}\sum_{k\in\{1,\dots,K\}\setminus I}(\mu\times\mu)(C_k\times C_k). 
 		\end{align*}
 		This implies
 		$$\sum_{k\in\{1,\dots,K\}\setminus I}(\mu\times\mu)(C_k\times C_k)\le \eta^2\gamma/2. $$
 		Since $\min_{1\le k\le K}(\mu\times\mu)(C_k\times C_k)=\eta^2$, we have shown that 
 		$$K-|I|\le \gamma/2<1,$$ which proves the claim, i.e., \eqref{claim4-1} holds.
 		
 	\smallskip
 	
 			\noindent{\bf (3)$\Rightarrow$(1):} Assume that $\mu$ is equicontinuous in the mean along $\{F_n\}_{n\in  \mathbb{N}}$. 
			By definition, given any $\tau>0$, there exists a compact
subset $Q$ of $X$ with $\mu(Q)>1-\tau$ such that $Q$ is equicontinuous in the
mean along $\{F_n\}_{n\in  \mathbb{N}}$. 
For any $\epsilon>0$, there exists $\delta>0$ such that
 $$x,y\in Q, d(x,y)<\delta \Rightarrow d_{m_G|_{F_n}}(x,y)<\epsilon, \forall n\in \mathbb{N}.$$ 
 Since $Q$ is compact, there exists a sequence of points $\{x_1,\dots,x_K\}\subseteq X$ such that $$Q\subseteq\cup_{i=1}^{k}B(x_k, \delta)\subseteq\cup_{i=1}^{K}B_{d_{m_G|_{F_n}}}(x_k, \epsilon).$$
 Thus $$\mu(\cup_{k=1}^{K}B_{d_{m_G|_{F_n}}}(x_k, \epsilon))>\mu(Q)>1-\epsilon,\ \forall n\in \mathbb{N}.$$
 This leads to $$S_{m_G|_{F_n}}(X,\mu,G,\epsilon)\le K,\ \forall n\in \mathbb{N}.$$ 
Then $\mu$ has bounded max-mean-complexity along the F{\o}lner sequence $\{F_n\}_{n\in\N}$. 
 By Theorem \ref{thm-1}, $\mu$ has discrete spectrum, which completes the proof.
	\end{proof}

\section*{Acknowledgement}
We would like to thank Prof. Wen Huang and Prof. Song Shao for valuable remarks and discussions.
This work is partially supported by  National Key R\&D Program of China (No. 2024YFA1013602, 2024YFA1013600), NNSF of China (No. 12031019, 12371197, 12426201, 12471188, 12090012, 12171175, 12090010). 

\end{document}